\documentclass[journal]{IEEEtran}
\usepackage{cite}
\usepackage{amsmath}

\usepackage{setspace}
\usepackage{bm}
\usepackage{amssymb}
\usepackage{latexsym}
\usepackage{makecell}
\usepackage{graphicx}
\usepackage{url}
\usepackage{tabularx}
\usepackage{multirow}
\usepackage{float}
\usepackage{slashbox}
\usepackage{algorithm}
\usepackage{algorithmic}

\begin{document}

\title{Strategic Bidding and Equilibria in Coupled \\ Gas and Electricity Markets}

\author{Cheng~Wang,~\IEEEmembership{Student Member,~IEEE}, Wei~Wei,~\IEEEmembership{Member,~IEEE}, Jianhui~Wang,~\IEEEmembership{Senior Member,~IEEE}, Feng~Liu,~\IEEEmembership{Member,~IEEE}, Shengwei~Mei,~\IEEEmembership{Fellow,~IEEE}
\thanks{This work is supported by the National Natural Science Foundation of China (51321005). The work of J. Wang is sponsored by the U. S. Department of Energy.(\emph{Correspond to: Feng Liu})} 
\thanks{C. Wang, W. Wei, F. Liu, and S. Mei are with the State Key Laboratory of Power Systems, Department of Electrical Engineering and Applied Electronic Technology, Tsinghua University, 100084 Beijing, China. (c-w12@mails.tsinghua.edu.cn, wei-wei04@mails.tsinghua.edu.cn, lfeng@mail.tsinghua.edu.cn, meishengwei@mail.tsinghua.edu.cn).}
\thanks{J. Wang is with Argonne National Laboratory, Argonne, IL 60439, USA (jianhui.wang@anl.gov).}}

\maketitle

\begin{abstract}
— The wide adoption of gas fired units and power-to-gas technology brings remarkable interdependency between natural gas and electricity infrastructures. This paper studies the equilibria of coupled gas and electricity energy markets driven by the strategic bidding behavior: each producer endeavours to maximizes its own profit subjecting to the market clearing outcome. The market equilibria is formulated as an equilibrium problem with equilibrium constraints (EPEC). A special diagonalization algorithm (DA) is devised, in which the unilateral equilibria of the gas or electricity market is found in the inner loop given the rival's strategy; the interaction of the two markets are tackled in the outer loop. Case studies on two test systems validates the proposed methodology.
\end{abstract}


\section*{Nomenclature}

\subsection{Sets and Indices}
\begin{IEEEdescription}[\IEEEusemathlabelsep\IEEEsetlabelwidth{$aaaaaa$}]
\item[$s\in S$]      Strategic electricity producer (SEP)
\item[$o\in O$]      Non-strategic electricity producer (Non-SEP)
\item[$i\in I_s$]    Units owned by SEP
\item[$j\in J_o$]    Units owned by non-SEP
\item[$l_p\in L_p$]      Power grid lines
\item[$n_p\in N$]      Power grid nodes
\item[$b\in B$]      Energy block of each generator
\item[$v\in V$]      Strategic gas producer (SGP)
\item[$w\in W$]      Non-strategic gas producer (Non-SGP)
\item[$m\in M_v$]    Gas wells owned by SGP
\item[$x\in X_w$]    Gas wells owned by non-SGP
\item[$c\in C$]      Active pipelines in gas network
\item[$l_g\in L_g$]      Passive pipelines in gas network
\item[$n_g\in N_g$]      Gas network nodes
\item[$z\in Z$]      Power-to-Gas facilities
\item[$d_p\in D_p$]      Electricity loads
\item[$d_g\in D_g$]      Gas loads
\item[$\bar{\varphi}_{n_g}^s$] Non-gas-fired units of SEP
\item[$\varphi_{n_g}^o$] Gas-fired units of non-SEP
\item[$\bar{\varphi}_{n_g}^o$] Non-gas-fired units of non-SEP
\item[$\varphi_{\{\cdot\}}(n_g)$]  Components connected to node $n_g$
\item[$\varphi_{\{\cdot\}_1}(n_g)$] Components whose head node is $n_g$
\item[$\varphi_{\{\cdot\}_2}(n_g)$] Components whose tail node is $n_g$
\item[$L_{l_p}(1,n_p)$] Power transmission lines whose head node is $n_p$
\item[$L_{l_p}(2,n_p)$] Power transmission lines whose tail node is $n_p$
\item[$\phi_{\{\cdot\}}(n_p)$]     Components connected to node $n_p$
\item[$\Theta(n_p)$]    Neighbouring nodes of node $n_p$
\item[$Line_{l_p}$]    Head and tail node set of line $l_p$
\item[$\phi^{-1}(\cdot)$] Connection node of components in power system
\item[$\varphi^{-1}(\cdot)$] Connection node of components in gas system
\item[$\varphi_1^{-1}(\cdot)$] Head node of components in gas system
\item[$\varphi_2^{-1}(\cdot)$] Tail node of components in gas system
\end{IEEEdescription}

\subsection{Parameters}
\begin{IEEEdescription}[\IEEEusemathlabelsep\IEEEsetlabelwidth{$aaaaa$}]
\item[$\lambda_{ib}^s$]  Marginal cost of non-gas-fired unit of SEP
\item[$\lambda_{jb}^o$]  Marginal cost of non-gas-fired unit of non-SEP
\item[$P_{d_p}$]           Electricity load demand
\item[$P_{ib}^{s,max}$]  Upper limit of block $b$ of unit of SEP
\item[$P_{jb}^{o,max}$]  Upper limit of block $b$ of unit of non-SEP
\item[$F_{l_p}$]              Electricity transmission line capacity
\item[$B_{n_{p1}n_{p2}}$]       Admittance of transmission line
\item[$q_x^u$]            Upper limit of gas well of non-SGP
\item[$q_m^u$]            Upper limit of gas well of SGP
\item[$q^{max}_c$]        Capacity of active gas pipeline
\item[$q^{max}_{l_g}$]        Capacity of passive gas pipeline
\item[$q_{d_g}$]           Gas load demand
\item[$\zeta_{m}^v$]     Marginal cost of gas well of SGP
\item[$\zeta_{x}^w$]     Marginal cost of gas well of non-SGP
\item[$\tau$]             Energy conversion constant
\item[$\eta_{ib},\eta_{jb}$] Efficiency of units
\item[$\eta_{z}$]         Efficiency of power-to-gas facilities
\item[$\alpha^{max}$]     Maximal bidding price of SEP
\item[$\delta^{max}$]     maximal bidding price of SGP
\end{IEEEdescription}

\subsection{Variables}
\begin{IEEEdescription}[\IEEEusemathlabelsep\IEEEsetlabelwidth{$aaaaa$}]
\item[$\beta_{n_p}$]       Local marginal electricity price
\item[$\alpha_{ib}^s$]   Bidding price of SEP
\item[$P_{ib}^s$]        Cleared output of unit of SEP
\item[$P_{jb}^o$]        Cleared output of unit of non-SEP
\item[$\theta_{n_p}$]      Phase angle of power nodes
\item[$\varrho_{n_g}$]      Local marginal gas price
\item[$\delta_{m}^v$]    Bidding price of SGP
\item[$q_{m}^v$]         Cleared output of gas well of SGP
\item[$q_{x}^w$]         Cleared output of gas well of non-SGP
\item[$q_{c}$]           Gas flow in active pipeline
\item[$q_{l_g}$]           Gas flow in passive pipeline
\item[$P_{z}$]           Demand of P2G facilities
\end{IEEEdescription}

\section{Introduction}

\IEEEPARstart{O}vER the past decade, there has been remarkable interest in the utilization of gas fired units due to its high efficiency and low carbon emission, thanks to the breakthroughs in turbine technologies, say, the open-cycle gas turbine (OCGT) and the combined-cycle gas turbine (CCGT), as well as the dramatic decline in natural gas price, owing to the shale rock revolution \cite{ShaleGas}. Meanwhile, emerging P2G technology allows using excessive electricity produce gases in a vast form, such as hydrogen and natural gas \cite{P2G_test}, which can be stored in existing tanks in liquid or compressed format without intensive investments on upgrading the energy storage equipment. Such technologies greatly enhances the operating flexibility of power systems with volatile renewable energy integration. The potential bi-directional energy flows imposes stronger interdependency between gas and electricity systems, which calls for interdisciplinary research on the system operation and design.

Along this line, extensive researches can be found on the topic of coordinated operation of power and gas systems. For example, the optimal power and gas flow is analyzed in \cite{OGPF}; a stochastic unit commitment is proposed in \cite{Alabdulwahab_SUC_Gas_Power}, which aims to mitigate the variability of wind generation; a robust scheduling model considering uncertain wind power and demand response is investigated in \cite{Linquan_interval_demand_wind}. To absorb excessive wind power, P2G technology is comprehensively discussed in \cite{Clegg_P2G}, which allows a bi-directional interchange of energy. From the marketing perspective, the interdependency mentioned above impacts both gas and electricity markets: on one hand, the price of natural gas will influence the cost of gas-fired units as well as gas demand in the electricity market; on the other hand, the price of electricity will influence the cost of P2G facilities as well as electricity demand in the gas market. Some inspiring work have been done to address the correlation between electricity and gas markets. The day-ahead strategic bidding behavior of one CCGT owner is discussed in \cite{Hantao_CCGT_MPEC}, taking operation constraints from both power and gas system into consideration, giving rise to a mathematical program with complementarity constraints (MPCC) formulation. The interaction of power and gas system models under a market environment in a medium-long scope is analyzed in \cite{Reneses_MidtermCo}, in which the single-level model includes power and gas operation constraints as well as the operational cost of either of the systems as its objective.

In the pool-based setting, the pool is cleared by a market operator (MO), an independent agency, given the bidding strategies of all participants and necessary data. The market power of strategic electricity producers has been well explored in existing studies since the pioneer work \cite{Hobbs-EM}. The optimal bidding strategy of a single producer is analyzed in \cite{Ruiz_MPEC}, resulting in a MPCC formulation. The optimal bidding strategies of multiple strategic producers are discussed in \cite{Ruiz_EPEC}, leading to an EPEC formulation, and is solved by a stationary point method. The strategic behavior of wind power producers are studied in \cite{QiaoWei_WindPool}, the resulting EPEC is solved by diagonalization method. Similar to the pool-based power market, the gas pool can be cleared by a gas market operator as well. In this paper, we consider several energy producers who trade their resources strategically, in coupled electricity and gas markets. The strategic producers seek to maximize their own profits by bidding their offering prices, and the market will reach an equilibrium.

As one of the first few attempt, this paper proposes a market mechanism that allows di-directional energy trading between independently cleared power and natural gas markets. An EPEC model is established to study the optimal bidding strategies and market powers of energy producers. Moreover, a dedicated diagonalization algorithm is designed to compute the market equilibrium, in which the inner level provides the equilibria associating with either the gas or electricity market while fixing the exogenous variable from the other market; the outer level coordinate the bidding strategies of both markets.

The rest of this paper is organized as follows. The basic settings and market mechanism are elaborated in Section~II. The EPEC formulation of the optimal bidding problem in coupled energy markets is presented in Section III. The  diagonalization algorithm is introduced in Section IV. To validate the proposed model and algorithm, several numerical results on two testing systems  are shown in Section V. Finally, Section VI draws the conclusion.

\section{Market Settings}
\subsection{Pool-based Market Mechanism}

\begin{figure}[H]
\centering
  \includegraphics[width=0.38\textwidth]{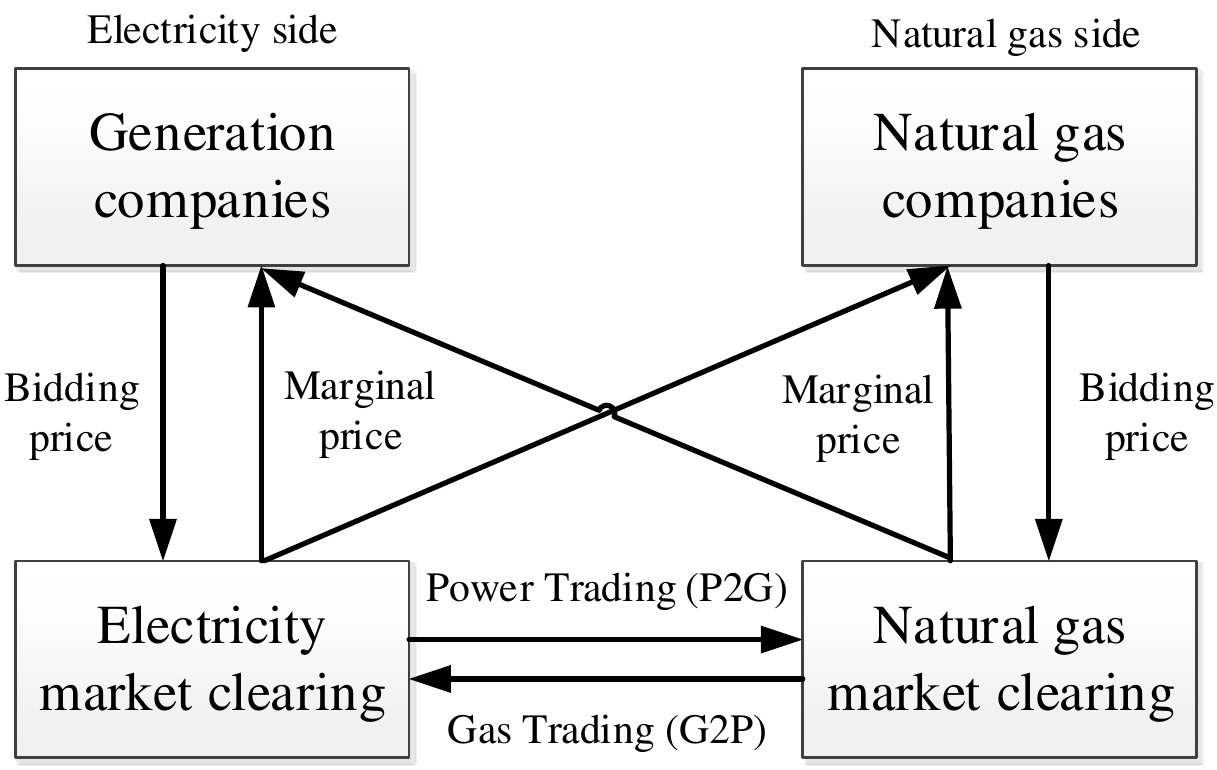}
  \caption{Electricity and natural gas market framework.}
  \label{fig:co_clearing}
\end{figure}

At the electricity side, SEPs such as generation companies receive local marginal gas price (LMGP) from the gas market operator (GMO), and then bid their offering price to the electricity market operator (EMO). The EMO clears the electricity market in terms of social welfare maximization, and the local marginal electricity price (LMEP) as well as the gas demand of each SEP become clear. At the gas side, strategic gas producers (SGPs) receive LMEP from the EMO, and then bid their offering price to the GMO. The GMO clears the gas market in terms of social welfare maximization, and the LMGP as well as electricity demand of each P2G facility become clear. It is apparent that the electricity market and the gas market impacts each other through bi-directional gas-electricity transitions. The schematic diagram of the coupled energy markets are shown in Fig. \ref{fig:co_clearing}.

\subsection{Assumptions and Simplifications}

The main assumptions made in the proposed model are summarized as follows.

1) General settings: i. the SEPs and SGPs are paid at LMEP and LMGP, respectively. Meanwhile, upper bounds are imposed on the bidding prices. ii. The demands in both markets are non-elastic. For detailed price-responsive load model, one can refer to \cite{Ruiz_EPEC}. iii. in the current formulation, we have adopted a static setting with a single period. The multi-period case incorporating inter-temporal constraints can be model in a similar way to that of \cite{Arroyo_DetailedMarket}.

2) For the power network: i. lossless DC power flow model is adopted. ii. Voltage magnitude at each bus is assumed to be close to 1.

3) For the gas network: i. we use a simplified steady-state model, in which the linepack is neglected. Moreover, we assume the gas flow in pipelines are directly controllable and omit the gas nodal pressure variable, leading to a linear approximation which is also used in \cite{Reneses_MidtermCo}, \cite{Gabriel_GasMarket}, \cite{Abrell_CoInvestment}. This simplification is reasonable when there are enough regulating devices. ii. We adopt the simplified compressor model in \cite{Wolf_BNG}, for detailed modeling of the compressor, please refer to \cite{Hantao_CCGT_MPEC}. iii. We don't incorporate gas storage devices in current formulation due to the static setting. For gas storage in a multi-period model, one can refer to \cite{Reneses_MidtermCo}. iv. We assume the P2G facilities to be non-strategic.

\section{Problem Formulation}
\subsection{EPEC model for electricity market}

In the electricity market, all SEPs aim to maximize their individual profit, leading to the following formulation

\begin{spacing}{0.8}
\begin{equation}
\begin{split}
\label{obj_mpec}
  \forall s\in S,~\min_{\alpha_{ib}^{s}}(&\sum_{(i\in \bar{\varphi}^{s}_u))b}\lambda_{ib}^{s}P_{ib}^{s}+\sum_{(i\in \varphi_i(n_g))b}\tau\varrho_{n_g}P_{ib}^{s}/\eta_{ib}\\&-\sum_{(i\in\phi_i(n_p))b}\beta_{n_p}P_{ib}^{s})
\end{split}
\end{equation}

\begin{equation}\label{strategic_offering}
~s.t.~~0\le \alpha_{ib}^{s} \le \alpha_{i(b+1)}^{s}\le \alpha^{max},~\forall i,\forall b<B
\end{equation}

\begin{equation}
\begin{split}
\label{obj_lower}
  P_{ib}^{s}, \beta_{n_g}\in \arg \{ &\min_{P_{ib}^s,P_{jb}^o,\theta}(\sum_{sib}\alpha_{ib}^{s}P_{ib}^{s}+\sum_{o(j\in\bar{\varphi}^o_{n_g})b}\lambda_{jb}^oP_{jb}^o\\
  &+\sum_{o(j\in\varphi_{n_g}^o)b}\tau\varrho_{n_g}P^{o}_{jb}/\eta_{jb})
\end{split}
\end{equation}

\begin{equation}\label{strategic_capacity}
~s.t.~~0\le P_{ib}^s\le P_{ib}^{s,max}: \beta_{ib}^{s,min},\beta_{ib}^{s,max}, \forall s,\forall i,\forall b
\end{equation}

\begin{equation}\label{non-strategic_capacity}
  0\le P_{jb}^o\le P_{jb}^{o,max}: \beta_{jb}^{o,min},\beta_{jb}^{o,max}, \forall o,\forall j,\forall b
\end{equation}

\begin{equation}
\begin{split}
\label{balance}
&~\sum_{s(i\in\phi_i(n_p))b}P_{ib}^s+\sum_{o(j\in\phi_j(n_p))b}P_{jb}^o-\sum_{z\in\phi_z(n_p)}P_{z}-\\ &\sum_{d\in\phi_d(n_p)}P_{d_p}=\sum_{k\in\Theta(n_p)}B_{n_pk_p}(\theta_{n_p}-\theta_{k_p}): \beta_{n_p},~\forall n_p
\end{split}
\end{equation}

\begin{equation}
\begin{split}
\label{transmission}
-F_{l_p}\le & B_{n_{p1}n_{p2}}(\theta_{n_{p1}}-\theta_{n_{p2}})\le F_{l_p}: \beta_{l_p}^{min},\beta_{l_p}^{max}\\ & n_{p1},n_{p2}\in Line_{l_p},\forall l_p
\end{split}
\end{equation}

\begin{equation}\label{Phase}
  -\pi \le \theta_{n_p} \le \pi: \varepsilon_{n_p}^{min},\varepsilon_{n_p}^{max}, \forall n_p
\end{equation}

\begin{equation}\label{Phase_1}
  \theta_{n_p}=0: \varepsilon_{1}, n_p=1~~\}
\end{equation}
\end{spacing}

\begin{spacing}{1.5}
\end{spacing}

In this formulation, the upper level bidding problem consists of (1)-(2) as well as the lower level electricity market clearing problem (3)-(9).  Objective function (\ref{obj_mpec}) represents the opposite profit of strategic producer $s$, where the first two components are the generation cost of non-gas-fired units and gas-fired units, respectively, and the last component is the utility of all generating units. (\ref{strategic_offering}) imposes upper and lower bounds on the offering price, as well as the monotonicity  with respect to the energy block. The LMEP $\beta_{n_p}$ and the energy contract $P^s_{ib}$ are determined by (\ref{obj_lower})-(\ref{Phase_1}), which represents the electricity market clearing process. Objective (\ref{obj_lower}) represents the opposite social welfare, where the first component is the generation expense of strategic producers and the last two component are the generation cost of non-gas-fired and gas-fired units owned by non-strategic producers, respectively. Constraints (\ref{strategic_capacity}) and (\ref{non-strategic_capacity}) are the energy block capacity of all type of units. Equality (\ref{balance}) depicts the nodal power balancing condition. Inequality (\ref{transmission}) is the flow limit of transmission lines. Inequality (\ref{Phase}) describes the upper and lower phase angle limits of the complex bus voltage. (\ref{Phase_1}) sets the phase angle for the reference node. Dual variables are indicated at the corresponding equations following a colon.

\subsection{EPEC model for gas market}

In the gas market, all SGPs aim to maximize their individual profit, leading to the following formulation

\begin{spacing}{0.8}
\begin{equation}\label{gas_obj_epec}
  \forall v\in V,~\min_{\delta_{m}^{v}}(\sum_{m}\zeta_{m}^{v}q_{m}^{v}-\sum_{(m\in\varphi_m(n_g))}\varrho_{n_g}q_{m}^{v})
\end{equation}

\begin{equation}\label{gas_con_upper}
  s.t.~~0\le \delta_{m}^{v}\le \delta^{max}, \forall m
\end{equation}

\begin{equation}
\begin{split}
\label{gas_obj_lower}
  q_{m}^{v},\varrho_{n_g}\in \arg \{&\min_{q_{m}^v,q_{x}^w,q_{c},q_{l_g},P_{z}}(\sum_{vm}\delta_{m}^{v}q_{m}^{v}+\sum_{wx}\zeta_{x}^wq_{x}^w+\\
  &\sum_{(z\in\phi_z(n_p))}\beta_{n_p}P_{z})
\end{split}
\end{equation}

\begin{equation}\label{gas_lower_con1}
  s.t.~~0\le q_{m}^v \le q_m^u:~\rho_{m}^{min},\rho_{m}^{max}~\forall m, \forall v
\end{equation}

\begin{equation}\label{gas_lower_con2}
  0\le q_{x}^w \le q_x^u:~\rho_{x}^{min},\rho_{x}^{max}~\forall x, \forall w
\end{equation}

\begin{equation}\label{gas_lower_con3}
  -q_{l_g}^{max}\le q_{l_g} \le q_{l_g}^{max}:~\rho_{l_g}^{min},\rho_{l_g}^{max}~\forall l_g
\end{equation}

\begin{equation}\label{gas_lower_con4}
  0\le q_{c} \le q_c^{max}:~\rho_{c}^{min},\rho_{c}^{max}~\forall c
\end{equation}

\begin{equation}\label{gas_lower_p2g}
  0\le P_{z}:~\rho_{z}~\forall z
\end{equation}

\begin{equation}
\begin{split}
\label{gas_lower_con5}
&\sum_{v(m\in\varphi_m(n_g))}q_{m}^v+\sum_{w(x\in\varphi_x(n_g))}q_{x}^w+\sum_{z\in\varphi_z(n_g)}\tau \eta_{z}P_{z}=\\
&\sum_{\{\cdot\}\in \varphi_{\{\cdot\}_1}(n_g)}q_{\{\cdot\}}-\sum_{\{\cdot\}\in\varphi_{\{\cdot\}_2}(n_g)}q_{\{\cdot\}}+\sum_{d_g\in\varphi_{d_g}(n_g)}q_{d_g}\\
&+\sum_{s(i\in\varphi_i(n_g))b}\tau P_{ib}^{s}/\eta_{ib}+\sum_{o(j\in\varphi_j(n_g))b}\tau P_{jb}^{o}/\eta_{jb}:\varrho_{n_g} \\
&~\{\cdot\}=\{c,l_g\},~\forall n_g~\}
\end{split}
\end{equation}
\end{spacing}

\begin{spacing}{2}
\end{spacing}

In this formulation, the upper level bidding problem consists of (10)-(11) as well as the lower level gas market clearing problem (12)-(18). The objective function (\ref{gas_obj_epec}) represents the opposite profit of strategic gas producer $v$, where the first component is the production cost and the second component is the utility of all the gas production. (\ref{gas_con_upper}) restricts the upper and lower bound of gas offering price. The LMGP $\varrho_{n_g}$ and the gas contract $q^v_m$ are determined by (\ref{gas_obj_lower})-(\ref{gas_lower_con5}), which represents the gas market clearing process. Objective function (\ref{gas_obj_lower}) defines the opposite social welfare, where the first two components are the generation expenses of strategic and non-strategic gas producers, respectively, and the last component is the electricity purchase cost of P2G facilities. Constraints (\ref{gas_lower_con1}) and (\ref{gas_lower_con2}) give the production limits of gas wells of strategic and non-strategic producers, respectively. In the natural gas network, pipelines without/with compressors are called passive/active. Gas flow capacities of passive and active pipelines are represented by (\ref{gas_lower_con3}) and (\ref{gas_lower_con4}), respectively. (\ref{gas_lower_p2g}) indicates a non-negative constraint on the electricity demand of P2G facilities. (\ref{gas_lower_con5}) is the nodal gas balancing condition. Dual variables are indicated at the corresponding equations following a colon.

\section{Solution Methodology}

In this section, we first present equivalent MPCC formulations for the bilevel bidding problems associated with each player in the upper level. Then a special diagonalization algorithm (DA) is introduced to solve the proposed EPEC models by iteratively solving a sequence of MPCCs, until a certain convergence criteria is met.

\subsection{Solution Methodology for MPCCs}

Take the electricity market as an example. For SEP $s^*$, if the LMGPs, the electricity demand from P2G facilities as well as the bidding strategies from other SEPs are fixed, the unilateral optimal bidding problem of SEP $s^*$ presented in (1)-(9) can be cast as the following MPCC

\begin{spacing}{0.5}
\begin{equation}
\begin{split}
\label{MPCC_SEP}
  \min_{\alpha_{ib}^{s^*}}(&\sum_{(i\in \bar{\varphi}^{s^*}_{n_g}))b}\lambda_{ib}^{s^*}P_{ib}^{s^*}+\sum_{(i\in \varphi_i(n_g))b}\tau\varrho_{n_g}P_{ib}^{s^*}/\eta_{ib}-\\&\sum_{(i\in\phi_i(n_p))b}\beta_{n_p}P_{ib}^{s^*})
\end{split}
\end{equation}

\begin{equation*}
~~s.t.~~(\ref{strategic_offering}), (\ref{strategic_capacity})-(\ref{Phase_1})
\end{equation*}

\begin{equation}
\label{Power_dual_1}
  \alpha_{ib}^s-\beta_{n_p}+\beta_{ib}^{s,max}-\beta_{ib}^{s,min}=0,~\forall i, \forall b, \forall s, n_p=\phi_i^{-1}
\end{equation}

\begin{equation}
\begin{split}
\label{Power_dual_4}
  &\lambda_{jb}^o-\beta_{n_p}+\beta_{jb}^{o,max}-\beta_{jb}^{o,min}=0,\\&~j\in \varphi_{n_g}^o, \forall b, \forall s, n_p=\phi_j^{-1}
\end{split}
\end{equation}

\begin{equation}
\begin{split}
\label{Power_dual_7}
  &\tau\varrho_{n_g}/\eta_{jb}-\beta_{n_p}+\beta_{jb}^{o,max}-\beta_{jb}^{o,min}=0,\\&~j\in \bar{\varphi}_{n_g}^o, \forall b, \forall s, n_p=\phi_j^{-1}
\end{split}
\end{equation}

\begin{equation}
\begin{split}
\label{Power_dual_10}
  &\sum_{k_p\in\Theta(n_p)}B_{nk}(\beta_{n_p}-\beta_{k_p})+\varepsilon_{n_p}^{max}-\varepsilon_{n_p}^{min}-(\varepsilon_{1})_{n_p=1}\\
  &+\sum_{n_p\in L_{l_p}(1,n_p)}B_{n_pk_p}(\beta_{l_p}^{max}-\beta_{l_p}^{min})-\\
  &\sum_{n_p\in L_{l_p}(1,n_p)}B_{n_pk_p}(\beta_{l_p}^{max}-\beta_{l_p}^{min})=0,~\forall n_p
\end{split}
\end{equation}

\begin{equation}\label{KKT_S_min}
  0\le P_{ib}^s\perp\beta_{ib}^{s,min}\ge 0,~\forall s, \forall i, \forall b
\end{equation}

\begin{equation}\label{KKT_S_max}
  0\le P_{ib}^{s,max}-P_{ib}^s\perp\beta_{ib}^{s,max}\ge 0,~\forall s, \forall i, \forall b
\end{equation}

\begin{equation}\label{KKT_O_min}
  0\le P_{jb}^o\perp\beta_{jb}^{o,min}\ge 0,~\forall o, \forall j, \forall b
\end{equation}

\begin{equation}\label{KKT_O_max}
  0\le P_{jb}^{o,max}-P_{jb}^o\perp\beta_{jb}^{o,max}\ge 0,~\forall o, \forall j, \forall b
\end{equation}

\begin{equation}
\begin{split}
\label{KKT_Trans_max}
  &0\le F_{l_p}-B_{n_{p1}n_{p2}}(\theta_{n_{p1}}-\theta_{n_{p2}})\perp \beta_{l_p}^{max}\ge0,\\ &~n_{p1},n_{p2}\in Line_{l_p}, \forall l_p
\end{split}
\end{equation}

\begin{equation}
\begin{split}
\label{KKT_Trans_min}
  &0\le B_{n_{p1}n_{p2}}(\theta_{n_{p1}}-\theta_{n_{p2}})+F_{l_p}\perp \beta_{l_p}^{min}\ge0,\\ &~n_{p1},n_{p2}\in Line_{l_p}, \forall l_p
\end{split}
\end{equation}

\begin{equation}\label{KKT_theta_max}
  0\le \pi-\theta_{n_p}\perp \varepsilon_{n_p}^{max} \ge0,~\forall n_p
\end{equation}

\begin{equation}\label{KKT_theta_min}
  0\le \theta_{n_p}+\pi\perp \varepsilon_{n_p}^{min} \ge0,~\forall n_p
\end{equation}
\end{spacing}

\begin{spacing}{1.5}
\end{spacing}
\noindent where, (\ref{Power_dual_1})-(\ref{KKT_theta_min}) are the KKT optimality conditions of the lower level electricity market clearing problem, the notation $0 \le a \perp b \ge 0$ represents the complementarity and slackness conditions $a \ge 0$, $b \ge 0$, and $ab=0$.

Similarly, for SGP $v^*$ in the gas market, if the LMEPs, the gas demand from gas-fired units as well as the bidding strategies from other SGPs are fixed, the unilateral optimal bidding problem of SGP $v^*$ presented in (10)-(18) can be cast as the following MPCC

\begin{spacing}{0.5}
\begin{equation}\label{SGP_MPCC}
  \min_{\delta_{m}^{v^*}}(\sum_{m}\zeta_{m}^{v^*}q_{m}^{v^*}-\sum_{(m\in\varphi_m(n_g))}\varrho_{n_g}q_{m}^{v^*})
\end{equation}

\begin{equation*}
~~s.t.~~(\ref{gas_con_upper}), (\ref{gas_lower_con1})-(\ref{gas_lower_con5})
\end{equation*}

\begin{equation}\label{gas_dual_1}
  \delta_{m}^v+\rho_{m}^{max}-\rho_{m}^{min}-\varrho_{n_g}=0,~\forall m, \forall v, n_g=\varphi^{-1}(m)
\end{equation}

\begin{equation}\label{gas_dual_2}
  \zeta_{x}^w+\rho_{x}^{max}-\rho_{x}^{min}-\varrho_{n_g}=0,~\forall x, \forall w, n_g=\varphi^{-1}(x)
\end{equation}

\begin{equation}
\begin{split}
\label{gas_dual_3}
  &\rho_{c}^{max}-\rho_{c}^{min}+\varrho_{n_{g1}}-\varrho_{n_{g2}}=0,\\&~\forall c, n_{g1}=\varphi^{-1}_1(c), n_{g2}=\varphi^{-1}_2(c)
\end{split}
\end{equation}

\begin{equation}
\begin{split}
\label{gas_dual_4}
  &\rho_{l_g}^{max}-\rho_{l_g}^{min}+\varrho_{n_{g1}}-\varrho_{n_{g2}}=0,\\&~\forall l_g, n_{g1}=\varphi^{-1}_1(l_g), n_{g2}=\varphi^{-1}_2(l_g)
\end{split}
\end{equation}

\begin{equation}\label{gas_dual_5}
  \beta_{n_p}-\rho_{z}-\tau\eta_z\varrho_{n_g}=0,\forall z, n_p=\phi^{-1}(z), n_g=\varphi^{-1}(z)
\end{equation}

\begin{equation}\label{gas_KKT_v_min}
  0\le q_{m}^v \perp \rho_{m}^{min} \ge0,~\forall m, \forall v
\end{equation}

\begin{equation}\label{gas_KKT_v_max}
  0\le q_m^u-q_{m}^v \perp \rho_{m}^{max} \ge0,~\forall m, \forall v
\end{equation}

\begin{equation}\label{gas_KKT_w_min}
  0\le q_{x}^w \perp \rho_{x}^{min} \ge0,~\forall x, \forall w
\end{equation}

\begin{equation}\label{gas_KKT_w_max}
  0\le q_x^u-q_{x}^w \perp \rho_{x}^{max} \ge0,~\forall x, \forall w
\end{equation}

\begin{equation}\label{gas_KKT_pas_min}
  0\le q_{l_g}+q_{l_g}^{max} \perp \rho_{l_g}^{min} \ge 0, \forall l_g
\end{equation}

\begin{equation}\label{gas_KKT_pas_max}
  0\le q_{l_g}^{max}-q_{l_g} \perp \rho_{l_g}^{max} \ge 0, \forall l_g
\end{equation}

\begin{equation}\label{gas_KKT_act_min}
  0\le q_{c} \perp \rho_{c}^{min} \ge 0, \forall c
\end{equation}

\begin{equation}\label{gas_KKT_act_max}
  0\le q_{c}^{max}-q_{c} \perp \rho_{c}^{max} \ge 0, \forall c
\end{equation}

\begin{equation}\label{gas_KKT_p2g}
  0\le P_{z} \perp \rho_{z} \ge 0, \forall z
\end{equation}
\end{spacing}

\begin{spacing}{1.5}
\end{spacing}

\noindent where, (\ref{gas_dual_1})-(\ref{gas_KKT_p2g}) are the KKT optimality conditions of the lower level gas market clearing problem.

There are two sorts of nonlinearities in the proposed MPCC: one is the summation of bilinear terms $\sum_{(i\in\phi_i(n_p))b}\beta_{n_p}P_{ib}^{s^*}$ in (\ref{MPCC_SEP}) and $\sum_{(m\in\varphi_m(n_g))}\gamma_{n_g}q_{m}^{v*}$ in (\ref{SGP_MPCC}), consisting of products of primal variable and dual variable; the other is the complementarity constraints, including (\ref{KKT_S_min})-(\ref{KKT_theta_min}) and (\ref{gas_KKT_v_min})-(\ref{gas_KKT_p2g}). The linearization methods are introduced below.

\subsubsection{Linearizing the primal-dual product}

linear expressions can be obtained to replace the bilinear terms in the original objective function by strong duality theorem and KKT complementarity constraints. The linearized objective functions are given below. Details can be found in the Appendix.
\begin{spacing}{0.5}
\begin{equation}
\begin{split}
\label{linear_power}
 &\min_{\alpha^{s^*}_{ib},P^{s^*}_{ib}}(\sum_{(i\in\varphi_i(n_g))b}\tau\varrho_{n_g}P_{ib}^{s^*}/\eta_{ib}+\sum_{(i\in\bar{\varphi}^{s^*}_{n_g}))b}\lambda_{ib}^{s^*}P_{ib}^{s^*}+\\
 &\sum_{n_p}\pi(\varepsilon_{n_p}^{max}+\varepsilon_{n_p}^{min})+\sum_{l_p}F_{l_p}(\beta_{l_p}^{max}+\beta_{l_p}^{min})\\
 &+\sum_{ib}(\sum_{s}P_{ib}^{s,max}\beta_{ib}^{s,max}+\sum_{o}P_{ib}^{o,max}\beta_{ib}^{o,max})\\
 &+\sum_{(s\neq s^*)ib}\alpha_{ib}^sP_{ib}^s+\sum_{o(j\in\bar{\varphi}^o_{n_g})b}\lambda_{jb}^oP_{jb}^o+\beta_{n_p}\sum_{d_p\in\phi_{d_p}(n_p)}P_{d_p}+\\
 &\sum_{o(j\in\varphi_{n_g}^o)b}\tau\varrho_{n_g}P^{o}_{jb}/\eta_{jb}-\beta_{ib}^{{s^*},max}P_{ib}^{{s^*},max}+\sum_{z\in\phi_z(n_p)}\beta_{n_p}P_{z})
\end{split}
\end{equation}

\begin{equation}
\begin{split}
\label{linear_gas}
&\min_{\delta_{m}^{v^*},q_{m}^{v^*}}(\sum_{m}\zeta_{m}^{v^*}q_{m}^{v^*}-\sum_{m}\rho_{m}^{max}q_{m}^{max}-\sum_{m}\delta_{m}^{v^*}q_{m}^{v^*}-\\
&\sum_{(v\neq v^*)m}\delta_{m}^{v}q_{m}^{v}-\sum_{(z\in\phi_z(n_p))}\beta_{n_p}P_{z}-\sum_{wx}\zeta_{x}^wq_{x}^w)-\\
&-\sum_{c}q_c^{max}\rho_{c}^{max}-\sum_{vm}\rho_{m}^{max}q_m^u-\sum_{wx}\rho_{x}^{max}q_x^u-\\
&\sum_{l_g}(\rho_{l_g}^{min}+\rho_{l_g}^{max})q_y^{max}+(\sum_{e\in\varphi_{e}(n_g)}q_{e}+\\
&\sum_{s(i\in\varphi_i(n_g))b}\tau P_{ib}^{s}/\eta_{ib}+\sum_{o(j\in\varphi_j(n_g))b}\tau P_{jb}^{o}/\eta_{jb})\varrho_{n_g}
\end{split}
\end{equation}
\end{spacing}

\subsubsection{Complementarity Constraint}

All complementarity constraints in KKT condition share the similar form of
\begin{equation}\label{general_KKT}
  0\le a-f \perp g \ge 0
\end{equation}
\noindent Where $a,f,g$ represent the primal variable, constant and dual variable, respectively. By introducing a binary variable $h$ and the following constraints, (\ref{general_KKT}) can be fully linearized.
\begin{equation}
\begin{split}
\label{linear_KKT}
  a-f&\le BigM\cdot h\\
  g&\le BigM(1-h)
\end{split}
\end{equation}
\noindent Where $BigM$ is a sufficient large positive number.

In view of the linearized objective functions and constraints, the MPCCs of SEP and SGP have been converted into mixed integer linear problems (MILPs), and are readily solvable by using commercial solvers, such as Cplex and Gurobi.

\subsection{The Nested Diagonalization Algorithm}

In this subsection, a nested diagonalization algorithm is proposed to find the equilibria in the coupled electricity-gas market. For notation brevity, we use $P_{s} (s\in S)$ and $G_{v} (v\in V)$ to represent the strategic bidding problem for SEPs and SGPs in MILP form, respectively. The aggregated problems $\{ P_s \}, \forall s$ and $\{ G_v\}, \forall v$ represent denote the EPECs of the electricity market and gas market, respectively.

\subsubsection{Inner Diagonalization Algorithms}

Two inner DAs are introduced to solve the EPEC models for electricity and gas market, respectively, the flow charts are given below.

\begin{algorithm}[h]
\caption{DA for EPEC of electricity market }
\label{alg:electricity}
\begin{algorithmic}[1]
\STATE Get current $P_{z},\varrho_{n_g}$ as input. Set $\alpha_{ib}^{s,0}=\alpha^{max}$ for all $s,i,b$. Set the maximum iteration number $r^{max}$ and convergence criterion $\epsilon$. $r=1,Flag_p=0$.
\STATE Do for $s^*=1,\dots,S$, regard $\alpha_{ib}^{s,0} (s\neq s^*)$ as the bidding price from others. Solve $P_{s^*}$ and obtain $\alpha_{ib}^{s^*,r}$. Let $\alpha_{ib}^{s^*,0}=\alpha_{ib}^{s^*,r}$.
\STATE If $|\alpha_{ib}^{s,r}-\alpha_{ib}^{s,r-1}|$$\le$$\epsilon\cdot\max\{\alpha_{ib}^{s,r},\alpha_{ib}^{s,r-1}\}$ for all $s$, then terminate, and report $\alpha_{ib}^{s,r}$. Else,
\STATE if $r=r^{max}$, set $Flag_p=1$, then quit and report the inner DA for electricity market fails to converge. Else, $r=r+1$ and return to Step 2.
\end{algorithmic}
\end{algorithm}

\begin{algorithm}[h]
\caption{DA for EPEC of gas market}
\label{alg:gas}
\begin{algorithmic}[1]
\STATE Get current $P_{ib}^s,P_{jb}^o,\beta_{n_p}$ as input. Set $\delta_{m}^{v,0}=\delta^{max}$ for all $v,m$. Set the maximum iteration number $r^{max}$ and convergence criterion $\epsilon$. $r=1,Flag_g=0$.
\STATE Do for $v^*=1,\dots,V$, regard $\delta_{m}^{v,0} (v\neq v^*)$ as the bidding price from others. Solve $G_{v^*}$ and obtain $\delta_{m}^{v^*,r}$. Let $\delta_{m}^{v^*,0}=\delta_{m}^{v^*,r}$.
\STATE If $|\delta_{m}^{v,r}-\delta_{m}^{v,r-1}|$$\le$$\epsilon\cdot\max\{\delta_{m}^{v,r},\delta_{m}^{v,r-1}\}$ for all $v$, then terminate, and report $\delta_{m}^{v,r}$. Else,
\STATE if $r=r^{max}$, set $Flag_g=1$, then quit and report the inner DA for gas market fails to converge. Else, $r=r+1$ and return to Step 2.
\end{algorithmic}
\end{algorithm}

\subsubsection{Outer Diagonalization Algorithm}

Similar to \cite{Jinghua_Power_heat}, we adopt an iterative-based algorithm to find the equilibria of the coupled energy markets, which is referred to as outer DA. The flow chart is given below.

\begin{algorithm}[h]
\caption{DA fot equilibria of the coupled markets}
\label{alg:electricity_gas}
\begin{algorithmic}[1]
\STATE Set $\varrho_{n_g}^{0}=\delta^{max},P_{z}^{0}=0$, $P_{ib}^{s,0}=0$,$P_{jb}^{o,0}=0$,$\beta_{n_p}^{0}=\alpha^{max}$. Set the maximum iteration number $r^{max}$ and convergence criterion $\epsilon$. $r=1$.
\STATE Call Algorithm 1 and obtain $\alpha_{ib}^{s}$. Solve problem (\ref{obj_lower})-(\ref{Phase_1}), and obtain $P_{ib}^{s,r},P_{jb}^{o,r},\beta_{n_p}^{r}$.
\STATE Call Algorithm 2 and obtain $\delta_{m}^{v}$. Solve problem (\ref{gas_obj_lower})-(\ref{gas_lower_con5}), and obtain $P_{z}^{r},\varrho_{n_g}^{r}$.
\STATE If $|P_{z}^{r}-P_{z}^{r-1}|$$\le$$\epsilon\cdot\max\{P_{z}^{r},P_{z}^{r-1}\}$ for all $z$, $|P_{ib}^{s,r}-P_{ib}^{s,r-1}|$$\le$$\epsilon\cdot\max\{P_{ib}^{s,r},P_{ib}^{s,r-1}\}$ for all $s$, $|P_{jb}^{o,r}-P_{jb}^{o,r-1}|$$\le$$\epsilon\cdot\max\{P_{jb}^{o,r},P_{jb}^{o,r-1}\}$ for all $o$, all holds, then quit, and report $\alpha_{ib}^{s},\delta_{m}^{v},P_z,P^s_{ib},P^o_{jb}$. Else,
\STATE If $r=r^{max}$, or $Flag_p=1$, or $Flag_g=1$, then quit and report the algorithm fails to converge. Else, $r=r+1$ and return to Step 2.
\end{algorithmic}
\end{algorithm}

Other than the proposed nested DA, traditional DA is also a straightforward method to find the equilibria in the gas-electric market. However, the computational cost of traditional DA is much higher than the nested DA, sometimes even fails to solve the EPEC model. The reason is that the updating variables in the outer DA are already equilibriums, say equilibriums in electricity market and gas market, respectively, while the updating variables
in traditional DA are not.

\section{Illustrative Example}
In this section, we present numerical experiments on two testing systems to show the effectiveness of the proposed model and algorithm. The experiments are performed on a laptop with Intel(R) Core(TM) 2 Duo 2.2 GHz CPU and 4 GB memory. The proposed algorithms are implemented on MATLAB with YALMIP toolbox. MILP is solved by Gurobi 6.5, the optimality gap is set as 0.1$\%$ without particular mention.

\subsection{6-Bus Power System with 7-Node Gas System}
\begin{figure}[!t]
\centering
  \includegraphics[width=0.40\textwidth]{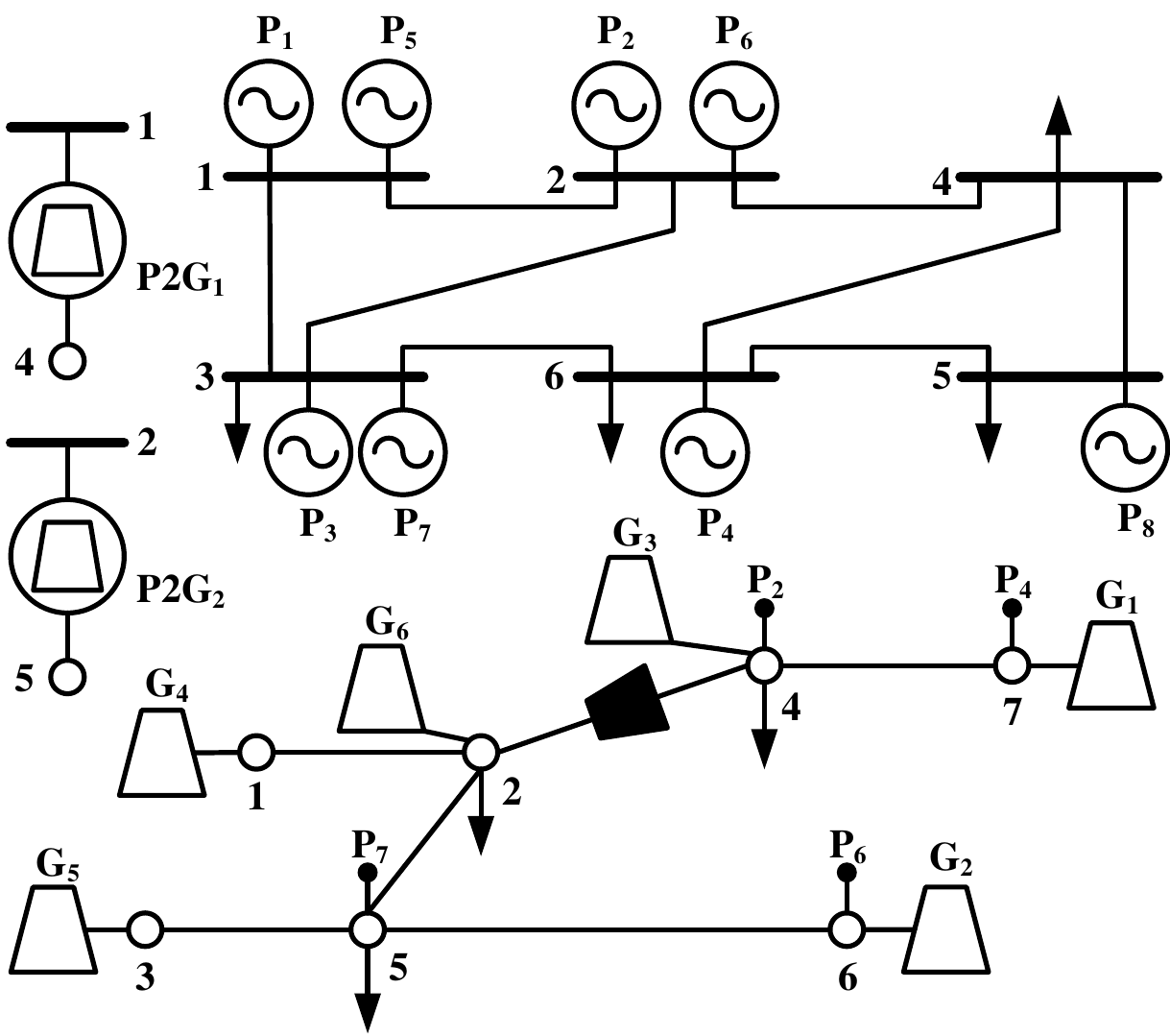}
  \caption{Topology of the hybrid test system.}
  \label{fig:power6gas7}
\end{figure}

Fig. \ref{fig:power6gas7} depicts the topology of the interconnected infrastructure. It has 4 gas-fired units ($P_2,P_4,P_6,P_7$), 4 traditional coal-fired units ($P_1,P_3,P_5,P_8$), 6 gas wells ($G_1-G_6$), 1 compressor, 2 P2G facilities ($P2G_1, P2G_2$), 4 power loads and 3 gas loads. Specifically, the head and tail node of the compressor is 4 and 2, respectively. The parameters of the system can be found in \cite{Power6Gas7}. The total electricity and gas demand are 525MW and 5520$Sm^3$, respectively. In the following cases, we have 4 strategic producers in the electricity market, who own generator $P_1$ to $P_4$, respectively, and the rest units are owned by a non-strategic producer. Also, we have 3 strategic producers in the gas market, who own gas well $G_1$ to $G_3$, respectively, and the rest gas wells are owned by a non-strategic producer. $\epsilon$ in the proposed diagonalization algorithm is $1\%$ and the maximum iteration number is 20.

\subsection{Uncongested Network}

The coupled energy market and its equilibria is firstly studied without considering congestions by neglecting capacity constraints, including those of power transmission lines and gas pipelines. The optimal bidding strategies and energy contracts are listed in Table~\ref{Tab:uncongested}, where both the superscript and the subscript are omitted for the price variables. The output of generators and gas wells which do not appear in Table I are 0. In the power network, all buses share the same LMEP, because the network is uncongested. The bidding price are equal to the corresponding LMEPs. The committed power from each block of generator can be obtained according to their capacity parameters. In the gas market, due to the existence of active pipeline as well as its radial topology, the LMGPs could be different at the head node and the tail nodes of an active pipeline. In this particular case, all the LMGPs are the same, and the strategic gas producers bid the same price equal to the LMGP. The energy exchange between electricity market and gas market are also listed in Table \ref{Tab:uncongested}. Intuitively, the power consumption of P2G facilities should be zero in an uncongested network as the equivalent gas price should be higher than the LMGPs in general, due to the dissipativity of energy conversion.
\begin{table}[!t]
\footnotesize
  \centering
  \newcommand{\tabincell}[2]{\begin{tabular}{@{}#1@{}}#2\end{tabular}}
  \caption{Results of uncongested network}\label{Tab:uncongested}
  \begin{tabular}{c|c|c|c|c|c|c}
  \hline
  \multicolumn{2}{c}{Power} & \multicolumn{2}{|c|}{Gas} & \multicolumn{3}{c}{Exchange}\\
  \hline
  \tabincell{c}{$\alpha$\\(\$/MWh)} & 16.67 & \tabincell{c}{$\delta$\\(\$/$Sm^3$)} & 0.9 & \multirow{8}{*}{\tabincell{c}{Gas\\($Sm^3/h$)}} & $P_2$ & 0\\
  \cline{1-4}
  \cline{6-7}
  \tabincell{c}{$\beta$\\(\$/MWh)} & 16.67 & \tabincell{c}{$\varrho$\\(\$/$Sm^3$)} & 0.9 & & $P_4$ & 2486\\
  \cline{1-4}
  \cline{6-7}
  \tabincell{c}{$P_4$\\(MW)} &143.5& \tabincell{c}{$G_1$\\($Sm^3/h$)} & 3300 & & $P_6$ & 2694\\
  \cline{1-4}
  \cline{6-7}
  \tabincell{c}{$P_6$\\(MW)} &157.6& \tabincell{c}{$G_2$\\($Sm^3/h$)} & 3000 & & $P_7$ & 1230\\
  \hline
  \tabincell{c}{$P_7$\\(MW)} &68.95 & \tabincell{c}{$G_3$\\($Sm^3/h$)} & 2630& \multirow{4}{*}{\tabincell{c}{Power\\(MW)}}&$P2G_1$& 0\\
  \cline{1-4}
  \cline{6-7}
  \tabincell{c}{$P_8$\\(MW)} &155 & \tabincell{c}{$G_4$\\($Sm^3/h$)} & 3000 & & $P2G_2$&0\\
  \hline
  \end{tabular}
\end{table}

\begin{figure}[!t]
\centering
  \includegraphics[width=0.40\textwidth]{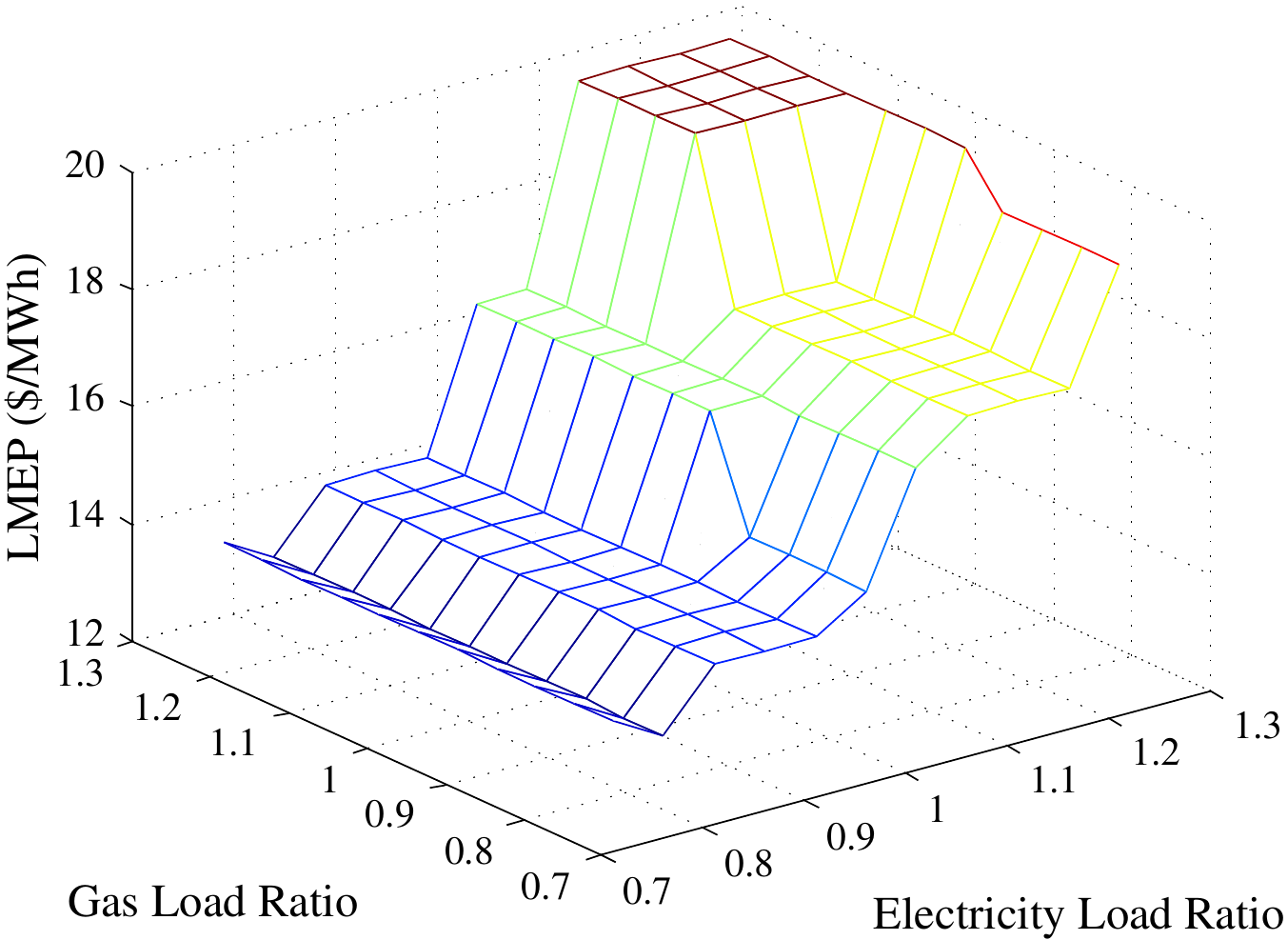}
  \caption{LMEPs under different load ratios.}
  \label{fig:Uncongested_Power}
\end{figure}

\begin{figure}[!t]
\centering
  \includegraphics[width=0.40\textwidth]{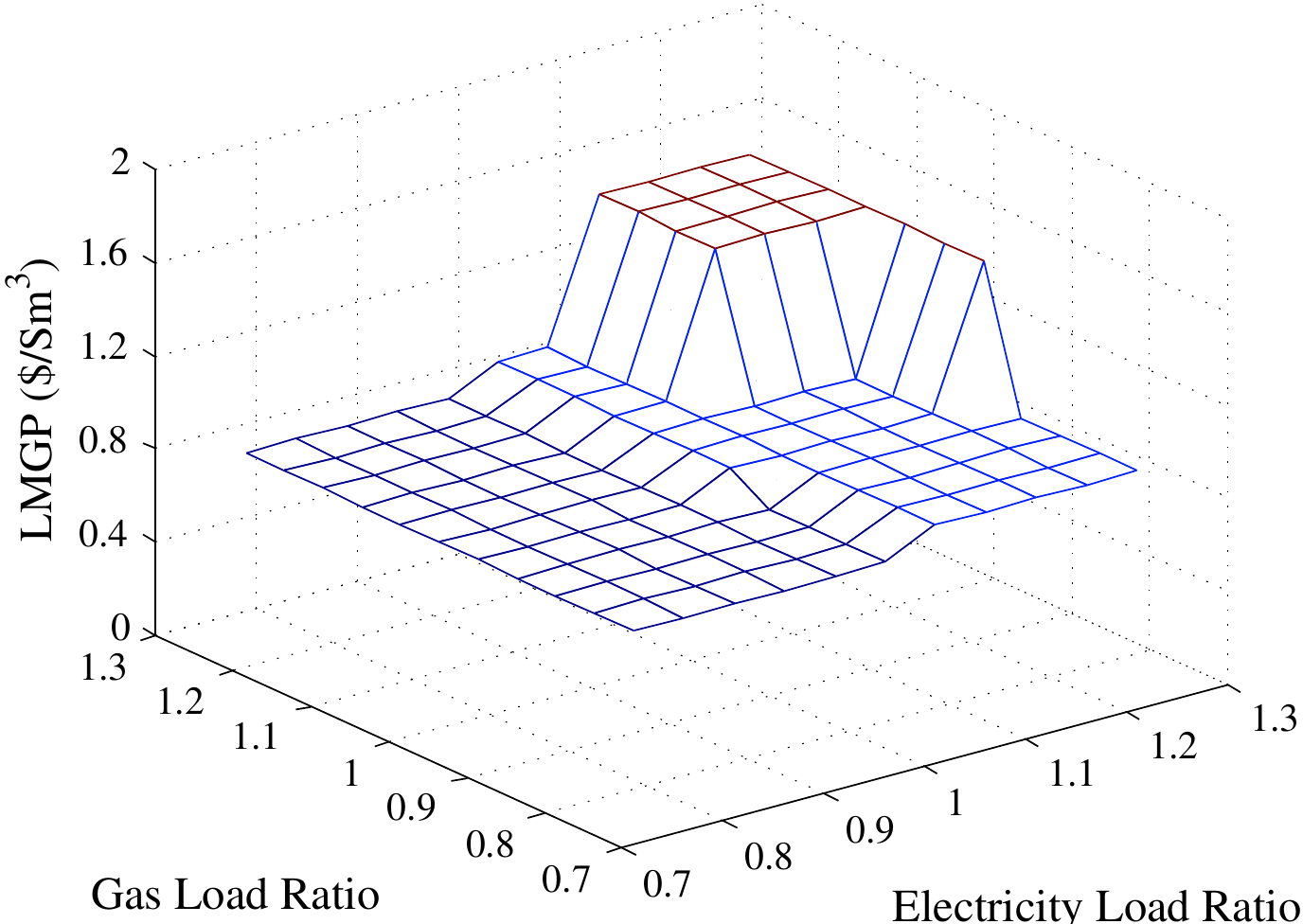}
  \caption{LMGPs under different load ratios. (Node:1-3, 5-6)}
  \label{fig:Uncongested_Gas_1}
\end{figure}

\begin{figure}[!t]
\centering
  \includegraphics[width=0.40\textwidth]{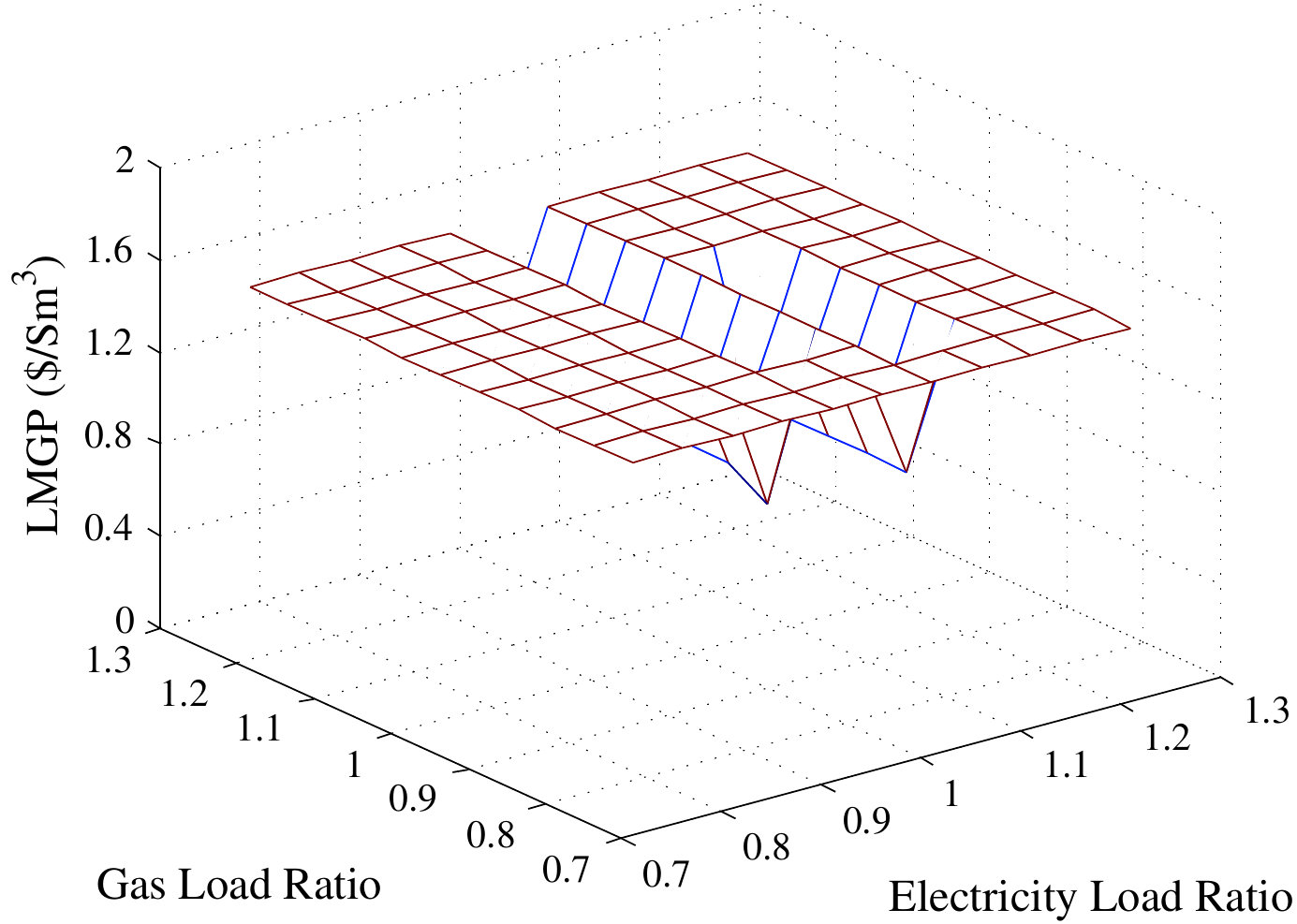}
  \caption{LMGPs under different load ratios.(Node:4,7)}
  \label{fig:Uncongested_Gas_2}
\end{figure}

Then we change the electricity and gas load by their load ratios, respectively. The corresponding nodal prices are shown in Fig. \ref{fig:Uncongested_Power}-Fig. \ref{fig:Uncongested_Gas_2}. From Fig. \ref{fig:Uncongested_Power}, it can be observed that the LMEP is always nondecreasing with the gas load ratio (GLR) when electricity load ratio is fixed. However, the relationship between LMEP and electricity load ratio (ELR) is quite complicated, whose general trend is increasing yet sometimes decreases, say there is a drop when ELR is 0.8, which is mainly caused by competition in the electricity pool as the LMGPs remain the same when ELR varies from 0.75 to 0.9. Further, from Fig. \ref{fig:Uncongested_Gas_1} and Fig. \ref{fig:Uncongested_Gas_2}, the LMGPs of two sides of active pipeline are different in most cases. Also, the general variation trend of Fig. \ref{fig:Uncongested_Power} and Fig. \ref{fig:Uncongested_Gas_1} are the same.

\subsection{Congested Network}

The analysis in the previous subsection is performed again with transmission capacity limits becoming mandatory, including those of power transmission lines and gas pipelines. The transmission capacity parameter can be found in \cite{Power6Gas7}. The results are listed in Table \ref{Tab:congested}. In the gas market, SGPs seek to bid higher prices due to congestion, which would give rise to the marginal cost increment of gas-fired units. In the electricity market, SEPs also seek to bid higher prices due to congestion, moreover, some of them are suffering from a higher marginal cost. Therefore, the bidding prices of SEP 4 as well as SGP 1,3 reach their maximum value, respectively, in the equilibria. In this regard, the MO of electricity market may try to dispatch more traditional generation assets to maximize social welfare. Thus, the total gas consumption from gas-fired units has a 2385 $Sm^3/h$ reduction compared with the uncongested case. Meanwhile, traditional units, say $P_1$ and $P_2$, are put into operation in this case while they are out of service in the uncongested case. It can be observed that the P2G facilities are not working as they are still not cost-effective according to the gas market clearing model.

From the maximum social welfare perspective, the production cost of the power network in the congested case is increased by \$221 compared with the uncongested case, from \$7,539 to \$7,760, and the production cost of the gas network in the congested case is increased by \$439 compared with the uncongested case, from \$9,537 to \$9,976.

\begin{table}[!t]
\footnotesize
  \centering
  \newcommand{\tabincell}[2]{\begin{tabular}{@{}#1@{}}#2\end{tabular}}
  \caption{Results of congested network}\label{Tab:congested}
  \begin{tabular}{c|c|c|c|c|c|c|c}
  \hline
  \multicolumn{4}{c}{Power} & \multicolumn{4}{|c}{Gas}\\
  \hline
  \multicolumn{2}{c}{$\beta$ (\$/MWh)} &\multicolumn{2}{|c}{Output (MW)}&\multicolumn{2}{|c|}{$\varrho$ (\$/$Sm^3$)}&\multicolumn{2}{c}{Output ($Sm^3/h$)}\\
  \hline
  1 & 15.42 & $P_1$ & 54.25 & 1 & 0.6 & 1 & 0\\
  \hline
  2 & 17.33 & $P_2$ & 55.15 & 2 & 0.9 & 2 & 534.4\\
  \hline
  3 & 13.51 & $P_3$ & 82.75 & 3 & 0.8 & 3 & 2251\\
  \hline
  4 & 36.18 & $P_4$ & 20.25 & 4 & 2 & 4 & 3000\\
  \hline
  5 & 38.09 & $P_5$ & 0 & 5 & 0.9 & 5 & 2000\\
  \hline
  6 & 40    & $P_6$ & 157.6 & 6 & 0.9 & 6 & 1760\\
  \hline
  \multicolumn{2}{c|}{$\alpha$ (\$/MWh)} & $P_7$ & 0 & 7 & 2 & \multicolumn{2}{c}{G2P ($Sm^3/h$)} \\
  \hline
  1 & 15.42 & $P_8$ & 155 & \multicolumn{2}{c|}{$\delta$ (\$/$Sm^3$)} & $P_2$ & 989\\
  \hline
  2 & 17.33 & \multicolumn{2}{c|}{P2G (MW)} & 1 & 2 & $P_4$ & 342 \\
  \hline
  3 & 13.51 & $P2G_1$ &0 & 2 & 0.9 & $P_6$ & 2694 \\
  \hline
  4 & 40 & $P2G_2$ & 0 & 3 & 2 & $P_7$ & 0\\
  \hline
  \end{tabular}
\end{table}

\subsection{118-Bus Power System with 20-Node Gas System}

In this subsection, we apply the proposed model and algorithm to a larger test system, which consists of a  modified IEEE 118-bus system and the Belgian high-calorific 20-node gas network. It possesses 30 gas-fired units, 24 traditional units, 10 gas wells, 3 compressors, 4 P2G facilities, 91 power loads and 9 gas loads. Refer to \cite{Power6Gas7} for the network topology and system data set. Table \ref{Tab:Time} summarizes the computational times with different numbers of strategic producers indicated by a tuple $(a,b)$, where $a$ is the number of SGPs and $b$ is the number of SEPs. The details of ownership of units and gas wells can be found in \cite{Power6Gas7}. From Table \ref{Tab:Time} we see that the computational time are always acceptable for the coupled energy market in this moderately sized system.
\begin{table}[!t]
\footnotesize
  \centering
  \newcommand{\tabincell}[2]{\begin{tabular}{@{}#1@{}}#2\end{tabular}}
  \caption{Computational time under different case settings }\label{Tab:Time}
  \begin{tabular}{c|c|c|c|c|c}
  \hline
  \multicolumn{6}{c}{Uncongested}\\
  \hline
  NSP & Time (s)& NSP & Time (s)& NSP & Time (s)\\
  \hline
  (3,3) & 99 & (4,3) & 117 & (5,3) & 48 \\
  \hline
  (3,4) & 167 & (4,4) & 195 & (5,4) & 203 \\
  \hline
  (3,5) & 259 & (4,5) & 209 & (5,5) & 271\\
  \hline
  \multicolumn{6}{c}{Congested}\\
  \hline
  NSP & Time (s)& NSP & Time (s)& NSP & Time (s)\\
  \hline
  (3,3) & 40 & (4,3) & 296 & (5,3) & 327\\
  \hline
  (3,4) & 39 & (4,4) & 313 & (5,4) & 324\\
  \hline
  (3,5) & 47 & (4,5) & 247 & (5,5) & 335\\
  \hline
  \end{tabular}
\end{table}

\section{Conclusion}

With the physical and economic interdependency between the natural gas and electricity systems becoming stronger, this paper propose an EPEC model that allows bi-directional gas-power trading, and a nested diagnalization algorithm to compute the market equilibria. Simulation results corroborate the effectiveness of the proposed method. The proposed methodology will be desired by policy makers and market designers.

\section*{Appendix}
For strategic producer $s^*$, multiplying both sides of (\ref{Power_dual_1}) by $P_{ib}^{s^*}$ and adding them together, yields following equation:
\begin{equation}
\begin{split}
\label{app_linear_1}
  \sum_{ib}(\alpha_{ib}^{s^*}+\beta_{ib}^{{s^*},max}-\beta_{ib}^{{s^*},min})P_{ib}^{s^*}-\sum_{(i\in\phi_i(n_p))b}\beta_{n_p}P_{ib}^{s^*}=0
\end{split}
\end{equation}

From (\ref{KKT_S_min})-(\ref{KKT_S_max}), we have
\begin{equation}
\begin{split}
\label{app_linear_2}
  &\beta_{ib}^{{s^*},min}P_{ib}^{s^*}=0 \Rightarrow \sum_{ib}\beta_{ib}^{{s^*},min}P_{ib}^{s^*}=0\\
  &\beta_{ib}^{{s^*},max}P_{ib}^{s^*}=\beta_{ib}^{{s^*},max}P_{ib}^{{s^*},max} \Rightarrow \\ &\sum_{ib}\beta_{ib}^{{s^*},max}P_{ib}^{s^*}=\sum_{ib}\beta_{ib}^{{s^*},max}P_{ib}^{{s^*},max}\\
\end{split}
\end{equation}

Substitute some terms in (\ref{app_linear_1}) with (\ref{app_linear_2}), we have
\begin{equation}
\begin{split}
\label{app_linear_3}
  \sum_{(i\in\phi_i(n_p))b}\beta_{n_p}P_{ib}^{s^*}=\sum_{ib}(\alpha_{ib}^{s^*}P_{ib}^{s^*}+\beta_{ib}^{{s^*},max}P_{ib}^{{s^*},max})
\end{split}
\end{equation}

From strong duality theorem, we have
\begin{equation}
\begin{split}
\label{app_strong_1}
&\sum_{ib}\alpha_{ib}^{s^*}P_{ib}^{s^*}+\sum_{(s\neq s^*)ib}\alpha_{ib}^sP_{ib}^s+\sum_{o(j\in\bar{\varphi}^o_{n_g})b}\lambda_{jb}^oP_{jb}^o+\\&
\sum_{o(j\in\varphi_{n_g}^o)b}\tau\varrho_{n_g}P^{o}_{jb}/\eta_{jb}=-\sum_{n_p}\pi(\varepsilon_{n_p}^{max}+\varepsilon_{n_p}^{min})\\
&-\sum_{ib}(\sum_{s}P_{ib}^{s,max}\beta_{ib}^{s,max}+\sum_{o}P_{ib}^{o,max}\beta_{ib}^{o,max})-\\
&\sum_{l_p}F_{l_p}(\beta_{l_p}^{max}+\beta_{l_p}^{min})+\sum_{oj}P_{j}^{o,min}\nu_{j}^{o,min}\\
&-\beta_{n_p}(\sum_{z\in\phi_z(n_p)}P_{z}+\sum_{d_p\in\phi_{d_p}(n_p)}P_{d_p})
\end{split}
\end{equation}

By Substituting (\ref{app_strong_1}) into (\ref{app_linear_3}), we arrive at the linear expression of $\sum_{(i\in\phi_i(n_p))b}\beta_{n_p}P_{ib}^{s^*}$ as shown in (47).

\vspace{3pt}
For strategic producer $v^*$, multiplying both sides of (\ref{gas_dual_1}) by $q_{m}^{v^*}$ and adding them together with respect to $m$, yielding the following equation:
\begin{equation}
\label{app_linear_2_1}
  \sum_{m}(\delta_{m}^{v^*}+\rho_{m}^{max}-\rho_{m}^{min})q_{m}^{v^*}-\sum_{(n_g=\varphi^{-1}(m))}\varrho_{n_g}q_{m}^{v^*}=0
\end{equation}

From (\ref{gas_KKT_v_min})-(\ref{gas_KKT_v_max}), we have (\ref{app_linear_2_2}). Substitute some terms in (\ref{app_linear_2_1}) with (\ref{app_linear_2_2}), we have (\ref{app_linear_2_3}). From strong duality theorem, we have (\ref{app_strong_2}). By substituting (\ref{app_strong_2}) into (\ref{app_linear_2_3}), we arrive at the linear expression of $\sum_{(m\in\varphi_m(n_g))}\gamma_{n_g}q_{m}^{v*}$ as shown in (48).

\begin{equation}
\begin{split}
\label{app_linear_2_2}
  & q_{m}^{v^*}\rho_{m}^{min}=0 \Rightarrow~\sum_{m}q_{m}^{v^*}\rho_{m}^{min}=0\\
  & q_{m}^{v^*}\rho_{m}^{max}=q_{m}^{max}\rho_{m}^{max} \Rightarrow~\sum_{m}q_{m}^{v^*}\rho_{m}^{max}=\sum_{m}q_{m}^{max}\rho_{m}^{max}
\end{split}
\end{equation}
\begin{equation}\label{app_linear_2_3}
  \sum_{(n_g=\varphi^{-1}(m))}\varrho_{n_g}q_{m}^{v^*}=\sum_{m}\delta_{m}^{v^*}q_{m}^{v^*}+\sum_{m}\rho_{m}^{max}q_{m}^{max}
\end{equation}

\begin{equation}
\begin{split}
\label{app_strong_2}
&\sum_{m}\delta_{m}^{v^*}q_{m}^{v^*}+\sum_{(v\neq v^*)m}\delta_{m}^{v}q_{m}^{v}+\sum_{(z\in\phi_z(n_p))}\beta_{n_p}P_{z}+\\
&\sum_{wx}\zeta_{x}^wq_{x}^w=-\sum_{c}q_c^{max}\rho_{c}^{max}-\sum_{vm}\rho_{m}^{max}q_m^u-\\
&\sum_{l_g}(\rho_{l_g}^{min}+\rho_{l_g}^{max})q_{l_g}^{max}-\sum_{wx}\rho_{x}^{max}q_x^u+(\sum_{d_g\in\varphi_{d_g}(n_g)}q_{d_g}\\
&+\sum_{s(i\in\varphi_i(n_g))b}\tau P_{ib}^{s}/\eta_{ib}+\sum_{o(j\in\varphi_j(n_g))b}\tau P_{jb}^{o}/\eta_{jb})\varrho_{n_g}
\end{split}
\end{equation}


\bibliographystyle{IEEEtran}
\bibliography{IEEEabrv,refs}

%




\end{document}